\input amstex
\documentstyle{amsppt}

\loadeufb
\loadeusb
\loadeufm
\loadeurb
\loadeusm

\magnification =\magstep 1
\refstyle{1}
\NoRunningHeads

\topmatter
\title On a conjecture of H. Wu 
\endtitle

\author Robert  Treger \endauthor
\address Princeton, NJ 08540  \endaddress
\email roberttreger117{\@}gmail.com \endemail
\endtopmatter

\document

In an email (January, 2015), Sullivan has asked the author  how much one can do for coverings of  compact complex manifolds. In this short note, we derive a conjecture of H. Wu from the theorem in \cite{4}, provided the fundamental group $\pi_1(X)$ is residually finite.
The author is grateful to Dennis Sullivan for the question.

\proclaim{Conjecture (H. Wu)} 
 The universal covering $U$ of a compact Kahler manifold $X$ of dimension $n$ with negative sectional curvature is a bounded domain in $\bold C^n$.
\endproclaim

\proclaim{Proposition} Let $U$ be the uiversal covering of a compact Kahler manifold $X$ of dimension $n$ with negative sectional curvature and
residually finite $\pi_1(X)$. Then $U$   is a bounded domain in $\bold C^n$.
\endproclaim

\demo{Proof} We consider $U$ with the  Kahler metric induced from $X$.

According to Wu \cite{5} (see also \cite{2, Problem B(ii), p.\;45}), $U$ is Stein hence  $\pi_1(X)$ is large, i.e., $U$ contains no proper submanifolds of positive dimension.

 According to Ballmann and Eberlein \cite{1}, $\pi_1(X)$ is nonamenable.
 
We will show $X$ is projective.  According to Gromov,  $X$ is Kahler hyperbolic \cite{3, (0.3.A), p.\;265}. The canonical bundle $\eusm K_X$ is quasiample \cite{3, (0.4.C), p.\;267; (3.2.B), (3.2.B$'$) p.\;287}. Thus $X$ is Moishezon and  Kahler hyperbolic hence $X$ is projective. The same follows from a theorem of Kodaira as follows. Since  the sectional curvature of $X$ is negative, the Ricci form of the volume form is negative. Hence $\eusm K_X$ is ample by Kodaira.

 Now, we can apply the theorem in  \cite{4, Introduction}.

\enddemo

\enddocument